\documentclass{article}
\usepackage{graphicx}
 \usepackage{mathptmx}
\usepackage{amsmath, amstext,amssymb,amsfonts}
\usepackage[english]{babel}
\usepackage{delarray}
\usepackage{mathptmx}
\usepackage{amsmath, amstext,amssymb,amsfonts}
\usepackage[english]{babel}
\usepackage{delarray}

\DeclareFontFamily{U}{mathx}{\hyphenchar\font45}
\DeclareFontShape{U}{mathx}{m}{n}{
      <5> <6> <7> <8> <9> <10>
      <10.95> <12> <14.4> <17.28> <20.74> <24.88>
      mathx10
      }{}
\DeclareSymbolFont{mathx}{U}{mathx}{m}{n}
\DeclareFontSubstitution{U}{mathx}{m}{n}
\DeclareMathAccent{\widecheck}{0}{mathx}{"71}
\DeclareMathAccent{\wideparen}{0}{mathx}{"75}

\numberwithin{equation}{section}

\newcommand{\R}{\mathbb R} \newcommand{\N}{{\mathbb N}}
\newcommand{\Co}{\mathbb C}
\newcommand{\M}{M({\mathbb R})}  \newcommand{\B}{B({\mathbb R})}  \newcommand{\Bl}{B_{loc}({\mathbb R})}
\newcommand{\A}{\cal{A}}    \newcommand{\pd}{\,positive-definite\, }
\newcommand{\U}{US(\mathbb{R})} \newcommand{\Hn}{PD(\mathbb{R})}

\date{}
\title{On  functional calculus for Hermitian elements of Banach algebras: the norm and spectral radius}
\author{\large{ Saulius  Norvidas}}
\date{\footnotesize Institute of Data Science and Digital Technologies, Vilnius University, \\ Akademijos str. 4, Vilnius LT-04812, Lithuania\\
 ({\rm{e-mail: norvidas{@}gmail.com}})}
\begin{document}

\maketitle
 {{ {\bf Abstract}}}
Let\ ${\cal A}$\ be a complex unital Banach algebra. An element\ $a \in {\cal A}$\ is said to be  Hermitian,
if\ $ \| \exp (ita) \| =1$\ for all\ $t\in\R$.\   In the case of the algebra  of bounded linear operators in a Hilbert space this Hermitian
 property agrees with the ordinary  selfadjointness. If\  $a \in {\cal A}$\ is  Hermitian, then\ $\|a\|=|a|$,\ where\ $|a|$\ denotes the spectral
radius of\ $a$.\  A function\ $F: \mathbb{R}\to \Co$\ is called the universal symbol if\ $ \| F(a)\|= |F(a)|$\  for each\ ${\cal A}$\
and all Hermitian\ $a\in {\cal A}$.\  We characterize  universal symbols  in terms of \pd functions.



\section{ Introduction }
{\large{
 Let\ ${\cal A}$\ be a complex unital Banach algebra. This means  that\ ${\cal A}$\ is  a  complex Banach space
 and  simultaneously an algebra with  unit\ $\mathbf{1}_{\cal A} $\ such that $ \|\, \mathbf{1}_{\cal A}\, \| =1 $,\ and
 \ $ \|ab \|\le \|a \|\, \| b\|$ \
 for all\ $a $, $b\in {\cal A} $.\ Given\ $a\in {\cal A}$,\ we denote by\ $ \sigma (a) $\ the {\it spectrum} of\ $a$,\  and by\
  $|a|=\max\{ |\lambda|:\ \lambda\in  \sigma (a)\}$\  the  {\it spectral radius} of\ $a$.\ The  {\it numerical range}\ $V(a)$\ \ is given by
  \[
 V(a)=\{ \varphi(a): \varphi\in {\cal A}^{\ast},\  \varphi(\mathbf{1}_{\cal A})=\|\varphi\|=1\},
  \]
where\  ${\cal A}^{\ast}$\  denotes the dual space  of  the Banach space\  ${\cal A}$.\  It is well  known (see Bonsall and Duncan 1971 [1])\ , that
\begin{equation}
\sigma(a)\subset V(a),\quad a\in {\cal A}.
\end{equation}
An element\ $a $\ is said to be  \textit {Hermitian}, if\ $V(a)\subset \R $.\ This condition
is fulfilled if and only if
\begin{equation}
 \| \exp (ita) \| =1 \quad {\text{ for all}}\ t\in\R
\end{equation}
( see [1]).\ It follows from\ (1.2)\ that for\ $C^{\ast}$-algebras
  the Hermitian property agrees with the ordinary definition of selfadjointness.
Other important examples of Hermitian elements we can obtain by considering the operator\ $D=-i\frac{d}{dx}$\ in various
 spaces. Indeed, in  the Banach space\ $B_{\varrho}^p$,\ $1\le p\le\infty$,\ $\varrho>0$,\ of entire functions\
$f$\ of exponential type $\le\varrho$\ such that\ $f\left|_{\R}\right.\in L^p(\R)$\ (see Levin 1996 [8]), the operator\
${\text{exp}}(itD)$\ defines a translation by\ $t\in\R$,\ and consequently\ $D$\ is Hermitian.

\ \ \ \ \ \  Let\ $H({\cal A})$\ denote the set of all Hermitian elements of\ ${\cal A}$.\ Then\ $H({\cal A})$\ forms a
real Banach space. Hermitian elements inherit   some properties of selfadjoint operators.  In particular,
\begin{equation}
\|a^n\|=|a^n|,\quad n\in \mathbb{N},
\end{equation}
(see Bonsall and Duncan 1973 [2]).\  On the other hand, there exist algebras\ ${\cal A}$\
and\ $a\in H({\cal A})$\ such that \ $a^2\not\in H({\cal A})$.\
 There are  known some other than\ $F(x)=x^n$\ functions\ $F$\  such that\ $ \| F(a)\|= |F(a)|$\ for Hermitian\
$a$\ (see [2],  Gorin 1988 [3],   K\"onig 1977 [7], and  Sinclair 1971 [11]).\ In the present paper we provide a description of  the class of all\ $F: \mathbb{R}\to \Co$\
such that\ $ \| F(a)\|= |F(a)|$\  for each complex unital Banach algebra\ ${\cal A}$\ and all\ $a\in  H({\cal A})$.\
Some results presented below have been announced in Norvidas 1986 [9];\  their proofs are available in preprint Norvidas 1999 [10].

\ \ \ \ \ \ The article is organized as follows. In the next section, we give some
necessary definitions and facts related to a functional calculus for Hermitian elements. The section ends by exact formulations of
the main results.  Section 3 contains proofs.

\ \ \ \ \ \ The author is grateful to Vidmantas Bentkus for comments which have greatly improved an earlier version of the manuscript.

 \section{  Preliminaries and main theorems}

\ \ \ \ \ Let\ $\M$\ be the Banach  algebra of all finite complex  Borel measures on\ $\R$\ equipped
with the total variation norm\ $\|\mu\|$,\ and convolution as multiplication.   We define the
Fourier-Stieltjes transform of\ $\mu\in\M$  as
\[
\hat{\mu}(\xi)=\int_{\R}e^{i\xi x} d\mu(x),\quad \xi\in\R.
\]
The set of all\ $\hat{\mu}$\ will be denoted by\ $\B$.\  Note that\ $\B$\ is a Banach algebra under the norm\
$\|\hat{\mu}\|=\|\mu\|$\ (so called   Fourier-Stieltjes algebra). The closed ideal in\ $\M$\  of all measures absolutely continuous
with respect to the Lebesgue measure can be identified with\ $L^1(\R)$.\  The closed ideal\ $A(\R)=\{\hat{\mu}\in\B: \
\mu\in L^1(\R)\}$\ is called the Fourier algebra.

\ \ \ \ \ \ By\ $\Bl$\ we denote the algebra of all complex-valued functions\ $F$\ on\
$\R$\ such that\ $F$\ belongs locally to\ $\B$.\ The example of\ $F(\xi)=\xi$\ shows that the inclusion\ $\B \subset  \Bl  $\ is strict.

\ \ \ \ \ Suppose\ $a\in H(\A)$.\ By\ (1.1)\ we have\ $\sigma(a)\subset V(a) \subset\R$.\ Therefore  we can apply to Hermitian elements  the extended
holomorphic functional calculus. Let\ $F\in\Bl$.\ If\ $F=\hat{\mu}$\ in a real
neighborhood of\ $\sigma(a)$,\ then one can set
\begin{equation}
F(a)=\int_{\R}{\text{exp}}(i t a)\, d\mu(t).
\end{equation}
 It is  straightforward that  (see\ Gorin 1988 [3]\ and\ K\"onig 1977 [7]):

\ \ \ \ (i)\ if\ $F_1=F_2$\ in a  real neighborhood of\ $\sigma(a)$,\ then\ $F_1(a)=F_2(a)$;

\ \ \ \ (ii)\ if\ $F$\ can be extended  to a function\ $h$,\ holomorphic in a\ $\mathbb{C}$-neighborhood of\ $\sigma(a)$,\ then\ $F(a)=h(a)$,
i.e.\ (2.1)\ extends the usual holomorphic functional calculus;

\ \ \  \ (iii)\ if\ $\sigma(F(a))$\ is the spectrum of the elements\ $F(a)\in\A$,\ then
\begin{equation}
\sigma(F(a))=F(\sigma(a)).
\end{equation}

 In view  of\ (2.1)\ and\ (i)--(iii),\ the elements of\ $\Bl$\ are called {\it symbols} for\ $H(\A)$.\ A symbol\ $F$\ is said to be {\it universal}, if
\begin{equation}
\|F(a)\|=|F(a)|
\end{equation}
for each unital Banach algebra\ $\A$\ and all\ $a\in H(\A)$.\ The set of all universal symbols we will denote by\
$US(\R)$.\ It follows from\ (1.2)\ and\ (1.3)\  that\ $F(x)={\text{exp}}(i\alpha x)$,\ $\alpha\in \R$,\ and\
$F(x)=x^n$,\ $x\in \mathbb{N}$,\ belong to\ $US(\R)$.\ The class\ $\U$\ is sufficiently  wide. For example,  if\
$F$\ is a polynomial with zeros lying on the imaginary axis, then\ $F\in\U$.\  However, so far no characterization of the class of universal symbols
are known. In this paper, we provide a complete description of\ $US(\R)$\
in terms of \pd functions. Let us introduce necessary  notation.

\ \ \ \ \ By\ $PD(\R)$\ we denote the class of all symbols\ $F$\  which admit   a positive-definite extensions from compact sets in the following sense.
By definition, \ $F\in PD(\R)$\ if for each compact set\ $K\subset \R$\ there exist\ $c\in\Co$,\ $x_0\in K$,\ and positive-definite function\ $\varphi$\
on\ $\R$\ such that
\begin{equation}
F(x)=c\varphi(x-x_0),\quad {\text{ for all}} \ x\in K.
\end{equation}
Note that we do not require the function\ $\varphi$\ in\ (2.4)\ to be continuous. Let\ $CH(\R)$\ stand for the class of\ $F$\  such that
the positive-definite function\ $\varphi$\  in\ (2.4)\  can be chosen to be continuous.  Obviously\ $CH(\R)\subset PD(\R)$.\
It is also known (see\ Gorin 1988 [3])\  that
\begin{equation}
US(\R)\subset PD(\R).
\end{equation}

\ \ \ \ \ \ {\bf Theorem 2.1.}\quad \textit{We have }
\begin{equation}
US(\R)=CH(\R)= PD(\R).
\end{equation}

\ \ \ \ \ \ Theorem 2.1  makes it possible to provide  new examples of universal symbols. For instance, if\ $F$\ is an even
polynomial with real coefficients with zeros satisfying\ $\pi/4\le |{\text{arg z}}|\le 3\pi/4$,\ then  \ $F\in\U$.\
 Also, by using the Polya theorem (for characteristic functions) it is possible to prove that, if\ $F$\ is a nonnegative
 and convex on\ $\R$,\ then\ $F\in\U$.\ However,  there is a number of  open  difficult problems related to the class\ $\U$.\
  For example,  it is still not known which complex polynomials of degree\ $n\ge 3$\ are universal. Theorem 2.1 clearly shows
  that any polynomial\ $F(x)=ax+b$,\ $a,b\in \Co$,\ is a universal symbol.  This fact was earlier established by Sinclair 1971 [11].\
  Gorin 2003 [4]\ described  universal symbols given by complex polynomials of the second degree.
  Excluding trivial cases, this problem   reduces to the polynomials
\[
F(x)=-x^2+2i x+t,\quad 0<t<2.
\]
Such polynomial is a universal symbol  if and only if \ $t\le t_0=(2\lambda^2_0- 1)/\lambda^2_0$,\
where\ $\lambda_0$\ is the root of\ $1-\lambda\cot \lambda=\lambda^2$\ in\ $(0; \pi)$.\ By the Lindeman
theorem,\ $t_0$\ is not algebraic and\ $t_0=1.867...$.\

\ \ \ \ \ \ In the following theorem we give some general properties of\ $F\in\U$.\
 For a symbol\ $F$\ put
\begin{equation}
m_F=\Bigl\{x\in\R: \ |F(x)|=\inf_{t\in\R}|F(t)|\Bigr\}.
\end{equation}

\ \ \ \ \ \ {\bf Theorem 2.2.}\quad \textit{Let\ $F\in\U$.\ The set\ $m_F$\ is connected. If \ $m_F\neq\R$\  and \ $U$\ is a connected
component of  \ $\R\setminus m_F$,\ then\ $|F|$\ is strictly monotone and unbounded on\ $U$. }

\ \ \ \ \ \ {\bf Corollary 2.3.}\quad \textit{A universal symbol\ $F\in\U$\ is bounded if and only if  there exist \ $\alpha\in\R$\ and\ $c\in\Co$\ such that }
\begin{equation}
F(x)=c e^{i\alpha x},\quad x\in\R.
\end{equation}

\ \ \ \ \ \ Theorems 2.1 and 2.2 allow us to describe some structural properties of\ $\U$.\  Indeed, it follows from
Theorem 2.1 and properties of positive-definite functions  that if\ $F\in\U$, then ${\bar F}$,\ $F^n$,\ $|F|^{2n}$,\
$n\in\N$,\ $cF$,\ $c\in\Co$,\ and\ $F(\alpha x+\beta)$,\ where\ $\alpha, \beta\in\R$,\ also are universal symbols. On
the other hand, if\ $F(x)=x^2+\alpha$,\ where\ $\alpha\in\R$,\ then\ $F$\ is a universal symbol if and only if\ $\alpha\ge 0$.\ By  Theorem 2.2,
we easily deduce that\ $\U$\ is not closed under  pointwise addition and multiplication of functions. In particular, the class\ $\U$\ is not a convex subset of\ $\Bl$.

 \section{  Lemmas and proofs of the theorems}

\ \ \ \ \ \ In what follows we  need to use the notion of a set of spectral synthesis and the concept of a positive-definite function;
here we give a brief account of some facts (for details see, for example,  Hewitt and Ross 1997 [6]).

\ \ \ \ \ \   Let\ $E$\ be a closed subset of\ $\R$\ and
denote by\ ${\cal I}(E)$\ the set of all\ $f\in L^1(\R)$\ such that\  ${\hat f}\left|_{E}\right.=0$.\ Let\ ${\cal J}(E)$\ be the closure of
the set of all\ $f\in L^1(\R)$\ such that\  ${\hat f}\left|_{W}\right.=0$\ for some open set\ $W$\ containing\ $E$\ (of course,\ $W$\
depends on the choice of\ $f$).\  Both\ ${\cal I}(E)$\ and\ ${\cal J}(E)$\ are closed ideals in\ $L^1(\R)$.\ If\  ${\cal I}(E)= {\cal J}(E)$,\
then\ $E$\ is said to be a {\it  spectral synthesis set} for the Fourier algebra\ $A(\R)$.\  The regular Banach algebra\
$A(\R)$\  satisfies  Ditkin's condition. Therefore,  if the boundary\ $\partial E$\  of\ $ E$\  contains no nonvoid perfect set, then\  $E$\ is
a spectral synthesis set  for\  $A(\R)$.\ The notion of a {\it  spectral synthesis set} for the Fourier-Stieltjes algebra\ $B(\R)$\
 is defined similarly.\  Recall that by Shilov-Wiener theorem, the algebra\ $A(\R)$\ coincides locally with\ $\M$.\
Therefore a compact subset\ $E$\ of\ $\R$\ is a  spectral synthesis set for\ $A(\R)$\ if an only if\ $E$\ is the same for\ $\M$.\

\ \ \ \ \ \   A function\ $\varphi:\R\to\Co$\ is said to be \pd if
\[
\sum_{j,k=1}^{n}\varphi(x_j-x_k)c_j{\bar c}_k\ge 0
\]
for all finite sets\ $x_1,\dots,x_n\in\R$\ and\ $c_1,\dots,c_n\in\Co$.\ Every \pd\ $\varphi$\ satisfies\ $\varphi(x)\le \varphi(0)$,\
$x\in\R$,\ and inequality
\begin{equation}
\bigl |\varphi (x_1)-\varphi (x_2)\bigr |^2 \le 2\varphi (0) \bigl|\varphi (0) - \Re (\varphi (x_1-x_2))\bigr|, \quad x_1, x_2\in\R.
\end{equation}
Continuous \pd\ $\varphi$\ is completely described by Bochner's theorem: such\ $\varphi$\ has exactly one representation\
$\varphi(x)={\hat \mu}(x)$,\ $x\in\R$,\ where\ $\mu$\  is a nonnegative measure in\ $\M$.\

\ \ \ \ \ \ {\bf Lemma 3.1}(see\ Hewitt and Ross 1997 [6], p.257).\quad \textit{Let\ $\varphi$\  be a positive-definite function on\ $\R$\ such that\
$\varphi(0)=1$\, and let\ $ X = \{x\in\R: |\varphi(x)| = 1\}$.\  Then\ $X$\  is a  subgroup of\ $\R$.\ For all\ $x\in X$\
and\ $y\in\R$,\  the equalities
\begin{equation}
\varphi(x+y)=\varphi(x)\varphi(y)
\end{equation}
hold. In particular,\ $\varphi$\  is a  character of the group\ $X$.}

\ \ \ \ \ \ Notice that if\ $\varphi$\  is a continuous \pd function, then\ $X$\ is a closed subgroup of\ $\R$.\

\ \ \ \ \ \ {\bf Remark.}\quad If a symbol\ $F$\  satisfies the conditions (2.4), then we say that\ $F$\
{\it admits a positive-definite extension from}\ $K$.\ In this case
\[
|F(x_0)|=\max_{x\in K}|F(x)|
\]
and we can takes\ $c=F(x_0)$.\  This notion  is well defined in the following sense:  if\  $x_1\in K $\ is another maxima point of\ $|F|$,\ and
\ $F\not\equiv 0$\ on\ $K$,\  then is easy to check (see  Lemma 3.1) that
\[
\psi(x):=\frac{F(x_0) \varphi(x-(x_0-x_1))}{F(x_1)}
\]
 is also  \pd function on\ $\R$\ and\ $F(x)=F(x_1)\psi(x-x_1)$\ for all\ $x\in K$.\ If\ $\varphi$\  in\ (2.4)\  can be chosen to be continuous, then
 we say that\ $F$\ {\it admits a characteristic extension from}\ $K$.

\ \ \ \ \ \ {\bf Lemma 3.2}.\quad \textit{Let\ ${\cal A}$\ be a complex unital Banach algebra,\ $a\in H({\cal A})$,\
and\ $F_1, F_2\in\Bl$.\ If there is a compact spectral synthesis set\ $K$\  for\ $B(\R)$\ such that\ $\sigma(a)\subset
K$\ and\ $F_1=F_2$\ on\ $K$,\  then\ $F_1(a)=F_2(a)$.}

\ \ \ \ \ \ {\bf Proof.}\quad If\ $F\in\Bl$\ and\ $F=\hat{\mu}$,\ $\mu\in\M$,\ in a real
neighborhood of\ $\sigma(a)$,\ then\ (2.1)\ implies
\begin{equation}
\|F(a)\|_{\cal A}=\|{\hat \mu}(a)\|_{\cal A}\le \|\mu\|_{\M}.
\end{equation}
Let\ $\nu_j$,\ $j=1, 2$,\ be measures in $\M$\ such that\ ${\hat \nu}_j=F_j$,\ $j=1, 2$,\ in some neighborhoods of\
$K$.\ Set\ $\nu=\nu_1-\nu_2$.\ Since\ $\sigma(a)\subset K$,\ we need only show that\ ${\hat \nu}(a)=0$.\ The condition\
$F_1\left|_{K}\right.=F_2\left|_{K}\right.$\ implies\ $\nu\in{\cal I}(K)$.\ Since\ $K$\ is a compact spectral synthesis
set for\ $B(\R)$,\ then  for every positive\ $\varepsilon$\ there is a measure\ $\eta\in {\cal J}(K)$\ such that\
$\|\nu-\eta\|<\varepsilon$.\ Thus\ ${\hat \eta}(a)=0$\ and\ (3.3)\ implies
\[
\|{\hat \nu}(a)\|=\|{\hat \nu}(a)-{\hat \eta}(a)\|\le\|\nu -\eta\|<\varepsilon.
\]
Therefore\ ${\hat \nu}(a)=0$,\ and so Lemma 3.2 is proved.

\ \ \ \ \ \ {\bf Lemma 3.3.}\quad \textit{Let}\  $F \in\Hn$,\  \textit{and suppose\ $| F|$\ assumes a relative maximum
value at a point\ $x_0\in\R$.\ Then\ $|F|$\ is a constant function in some neighborhood of\ $x_0$.
}

\ \ \ \ \ \ {\bf Proof.}\quad Set\  $G(x)=|F(x)|^2$.\ It is obvious that\ $G\in\Hn$\ and\ $x_0$\ is a maximum point of\
$G$.\ Suppose first that\ $x_0$\ is not a point of strongly   maximum of\ $G$.\ Since the case\ $G(x_0)=0$\ is obvious,
then we may, without loss of generality, assume that\ $G(x_0)=1$.\ Now there   exist\ $\omega>0$ and sequence\ $
\{y_m\}_{m=1}^{\infty}$\  of distinct elements of\ $[x_0-\omega, x_0+\omega]$\ such that\ $G(x)\le 1$\ for all\
$x\in [x_0-\omega, x_0+\omega]$,\ $G(y_m)=1$,\ and\ $\lim y_m=x_0$.\  By definition, put
\[
\varphi(x)=G(x+x_0),\quad x\in [-\omega, \omega].
\]
Since\ $G\in\Hn$, then\ $\varphi$\ can be continued to a \pd  function on\ $\R$.\ Therefore, the equalities\
$\varphi(y_m-x_0)=1$,\ $m=1,\dots,$\ imply by Lemma 3.1, that\ $\varphi\equiv 1$\     on a everywhere  dense subgroup
of\ $\R$.\ On the other hand, the continuity of\ $G\in\Bl$\
 on\ $\R$\ gives that\ $\varphi$\ is continuous on\ $[-\omega,\omega]$.\  Thus by\ (3.1),\ $\varphi$\ is continuous
 on the whole\ $\R$.\ Finally, we
 have that\ $\varphi\equiv 1$\ on\ $\R$,\ and consequently\  $G(x)= 1$,\ $x\in[x_0-\omega, x_0+\omega]$.

\ \ \ \ \ \ Let us prove that\  $G=|F|^2$\   cannot have the points of strongly relative maximum. Assume on the
contrary that there is such a point\ $x_0$.\ Because\ $\Hn$\ is   invariant under real translations, we can suppose
without restriction that\ $ x_0= 0$\ and\ $G(x_0) = 1$.\ Let\ $\omega$\ be any positive numbers such that
\[
G(x)<1,\quad x\in[-\omega,\omega]\setminus\{0\}.
\]
Recall that a point\ $t\in [-\omega,\omega]$\  is called right-visible for\ $G$\  on\ $[-\omega,\omega]$\ if\ $G(t)\ge
G(x)$\ for each\ $x\in[t,\omega]$.\ Let\ $\Sigma(G)$\ be the set of all right-visible points of\ $G$\ on\
$[-\omega,\omega]$.\  By the F. Riesz theorem\ $\Sigma(G)$\ is a closed subset of\ $[-\omega,\omega]$,\  and\ $G$\   is
monotone nonincreasing function on\ $\Sigma(G)$.\ Since\ $\Sigma\subset [0, \omega]$\ and\ $0, \omega\in \Sigma(G)$,\
 that only the following  cases are possible:\quad (i)\ There is\ $\omega_1\in(0,
\omega]$\ such that\ $[0,\omega_1]\subset \Sigma(G)$;\quad (ii)\ $\Sigma(G)$\ contains a sequence of intervals\
$\{[a_m,b_m]\}_1^{\infty}$\ such that\ $0<a_{m+1}<b_{m+1}<a_m<b_m\le\omega$\ and\ $\lim_{m\to \infty}b_m=0$.

\ \ \ \ \ \ Case (i). \  There exist\ $\alpha\in [-\omega,0)$\ and\ $\beta\in(0, \omega_1]$\ such that\
$0<2\beta-\alpha<\omega_1$\ and
\begin{equation}
 G(\alpha)=G(\beta)>0.
 \end{equation}
Put\ $D=\{\alpha\}\cup[\beta,\omega_1]$.\ For\ $x$\ belonging to\ $-\beta+D$,\ we define\
$\psi(x)=G(x+\beta)/G(\beta)$.\ Since\ $G\in\Hn$, then\ $\psi$\ can be continued to a \pd  function on\ $\R$.\ Now\
(3.4)\ implies\ $\psi(\alpha-\beta)=\psi(0)=1$.\  As\ $\psi$\ is \pd, we have\ $\psi(\beta-\alpha)=
{\overline{\psi(\alpha-\beta)}}=1$.\ Hence
\[
G(2\beta-\alpha)=G(\alpha)=G(\beta)>0.
\]
Because\ $G$\  is a monotone nonincreasing function on\ $\Sigma(G)$,\ it follows that\ $G$\ is a constant function on\
$[\beta, 2\beta-\alpha]$.\ From hence and from Lemma 3.1 we can deduce that the \pd extension of\ $\psi(x)$,\ $x\in[0,
\omega-\beta]$\ on\ $\R$\ is a constant function. This means that\ $G$\ is also a constant function on\
$[\beta,\omega_1]$.\  On the other hand, since\ $\beta$\ is an arbitrary point of\ $(0, \omega_1]$,\ we must have that\
$G$\ has not a strongly maximum  at\ $x_0=0$, a contradiction.

\ \ \ \ \ \ Case (ii). \  In this case\  \ $\Sigma(G)$\ contains a sequence of intervals\ $\{[a_m,b_m]\}_1^{\infty}$\
such that\ $0<a_{m+1}<b_{m+1}<a_m<b_m\le\omega$\ and\ $\lim_{m\to \infty}b_m=0$.\ By definition  of\ $\Sigma(G)$\ and
by the F. Riesz theorem mentioned above, we have\ $G(x)<G(a_m)$,\ $x\in(b_{m+1}, a_m)$,\ and
\begin{equation}
 G(a_m)=G(b_{m+1}), \quad m\in\N.
 \end{equation}
 If\ $x\in [b_{m+1}, b_m]$,\ then similarly as in the Case (i) we deduce that\ $G$\ is  a constant function
 on\ $[a_m, b_m]$.\  Thus\ (3.5)\ implies
 that\ $G$\ is the same constant function  on the whole set\ $\cup_{m=1}^{\infty}[a_m,b_m]$.\ Since\ $\lim_{m\to \infty}b_m=0$,\ we must have
 that\ $G$\ has not a strongly maximum  at\ $x_0=0$, a contradiction. This completes the proof of Lemma 3.3.

\ \ \ \ \ \ Recall that the quantity\ $m_F$ \  was defined in\ (2.7).

\ \ \ \ \ \ {\bf Lemma 3.4.}\quad \textit{Let}\  $F \in\Hn$.\  \textit{The set\ $m_F$\ is connected. If \ $m_F\neq\R$\
and \ $U$\ is a connected component of\ $\R\setminus m_F$,\ then\ $|F|$\ is strictly monotone and unbounded on\ $U$. }

\ \ \ \ \ \ {\bf Proof.}\quad The proof consists of three parts.

\ \ \ \ \ \ Part 1.\  Suppose\ $|F|$\ has a point\ $x_0$\ of a relative maximum. Let\ $G=|F|^2$.\ According to Lemma
3.3, there exists a neighborhood\ $\Gamma$\  of\ $x_0$\ such that\ $G$\ is a constant function on\ $\Gamma$.\ By\
$\Delta$\ denote the largest closed connected\ $\Gamma$\ of such type. The case\ $\Delta=\R$\ is trivial. Suppose without
restriction that\ $\Delta$\ is a bounded above, and\ $b$\ is the largest number such that\ $b\in \Delta$.\ Let us
prove that\ $G$\ is a strictly  monotone increasing function on\ $(b,\infty)$.\ To this end, we show at first that
\begin{equation}
G(t)>G(b)
 \end{equation}
 for all\ $t\in(b,\infty)$.\   Assume on the contrary that there is an $\lambda\in (b,\infty)$\ for which\ $G(b)\ge G(\lambda)$.\
 Then the function\ $\varphi(x):=G(x+b)$\  can be continued from the set\ $-b+\bigl([x_0, b]\cup \{\lambda\}\bigr)$\
  to a \pd  function on\ $\R$.\
 By Lemma 3.1, we obtain that\ $\varphi$\    is  a constant function on\ $\R$.\ Thus\ $G(\lambda)=G(b)$.\  Similarly
 we see that\
 $G(b)\le G(\lambda)$\ for all\ $t\in(b,\lambda]$.\ By Lemma 3.3, we have\ $G(t)=G(b)$\ for all\ $t\in (b, \lambda]$.\
 Therefore\ $\lambda\in\Delta$.\ This contradicts the fact that\ $b$\ is the largest
number such that\ $b\in \Delta$.\ Hence\ (3.6)\ holds.  Now we shall prove that\ $G$\ strictly monotone increase on\
$(b,\infty)$.\ Suppose on the contrary that there are\ $y_1, y_2\in (b,\infty)$\ such that\ $y_1<y_2$\ and
\[
G(b)<G(y_1)=\max_{t\in [b,y_2]}G(t).
\]
By Lemma 3.3 this implies that\ $G$\ is a constant function in some neighborhood of\ $y_1$.\ Thus\ $G$\  can be
continued from\ $[b, y_2]$\
 to a \pd  function on\ $\R$\ only as a constant function. Hence\ $G(b)=G(y_1)=G(y_2)$,\ contradicting to \ (3.6).\ This proves
 that\ $G$\ is a strictly monotone increasing function on\ $(b,\infty)$.\  In the case, if\ $\Delta$\ is also left-bounded and\ $a$\ is the  least
number such that\ $a\in \Delta$,\ then similarly we can show that\ $G$\ is a strictly monotone decreasing function on\
$(-\infty, a)$.\ Then\ $m_F=\Delta$.\
 This implies that\ $m_F$\ is a closed connected set.

\ \ \ \ \ \ Part 2.\  Suppose\ $|F|$\ has no points of a relative maximum. Then:  (i) either\ $|F|$\ has no points of a
relative minimum or (ii) \ $|F|$\ has a unique such point. If (i), then\ $m_{F}=\emptyset$\ and\ $|F|$\ is a strictly
monotone function  on the whole\ $\R$.\ If (ii), then\ $m_{F}=\{a\}$\,   and\ $|F|$\    strictly monotone decrease on\
$(-\infty, a)$,\  and strictly monotone increase on\ $(a,\infty)$.

\ \ \ \ \ \ Part 3.\   Let us prove that\ $|F|$\ is  unbounded on every connected component of\ $\R\setminus m_F$.\ As
was proved above, the set  \ $\R\setminus m_F$\ has at most two connected component. Let\ $U$\ be any such component.
Since the proof is essential the same in all cases, we may, without loss of generality, assume that\  $U=(a, \infty)$,\
$a>-\infty$.\ Assume on the contrary that \ $G=|F|^2$\ is bounded on\ $U$.\ Since\ $G$\ is a strictly  monotone
increasing function on\ $U=(a,\infty)$,\ then
\begin{equation}
\sup_{t\in [a,\infty)}G(t)=\lim_{t\to\infty}G(t)=\theta<\infty.
 \end{equation}
 Now fix a point\ $b\in(a,\infty)$,\  and  consider\ $G$\ on an arbitrary interval\ $[a, \sigma]$\
 such that\ $b\in(a, \sigma)$.\
Put, by definition\  $\varphi(x)=G(x+\sigma)$,\  $x\in [-\sigma+a, 0]$.\ Because\ $G$\ increase on\ $(a,\infty)$\ and
$G\in \Hn$\ it follows that\  $\varphi$\
 can be continued from\  $ [-\sigma+a, 0]$\    to a \pd  function on\ $\R$.\ Hence\ (3.1)\ and\ (3.7)\ imply
\begin{gather}
|G(a)-G(b)|^2=
\nonumber \\
|\varphi(-\sigma+a)-\varphi(-\sigma+b)|^2\le 2\theta|\varphi(0)-\varphi(a-b)|= 2\theta|G(\sigma)-G(\sigma+a-b)|.
\end{gather}
Now if in\ (3.7)\  $\sigma$\  turn  to infinite, then we have\ $G(a)=G(b)$.\ This contradicts the fact that\ $G$\ is a
strictly monotone increasing function on\ $(a,\infty)$,\  and so\ $G$\ is unbounded on\ $U=(a,\infty)$.

\ \ \ \ \ \ {\bf Proofs of Theorem 2.2 and Corollary 2.3.}\quad The statement of Theorem 2.2  follows now at once from
relation\ (2.5)\ and Lemma 3.4. Therefore, if  an universal symbol\ $F$\ is bounded on\ $\R$,\  then\ $m_F=\R$\ and\
$|F|$\ is a constant function on\ $\R$.\ If\ $|F|\not\equiv 0$,\ then we may, without loss of generality, suppose that\
$F(0)=1$.\ Since\ $F$\ admits a positive-definite extension from each interval\ $(-a,a)$, $a>0$,\ then Lemma 3.1
implies that
\begin{equation}
F(x+y)=F(x)F(y),\quad {\text{for all}}\  x,y\in \R.
 \end{equation}
Every\ $F\in\U\subset\Bl$\ is a continuous function on\ $\R$.\ Thus\ (3.9)\ implies that\ $F$\ is a continuous
characters on\ $\R$.\ This shows\ (2.8).

\ \ \ \ \ \ {\bf Proof of Theorem 2.1}\quad By\ (2.5)\ and trivial relation\ $CH(\R)\subset\Hn$,\ it suffices to show
that every\ $F\in\Hn$\ belongs to the classes\ $CH(\R)$\ and\ $\U$.\

\ \ \ \ \ \ Part 1. \ Suppose\ $F\in\Hn$\ and\ $K$\ is a compact subset of\ $\R$.\ Let us prove that\  $F$\ admits a
characteristic extension from\ $K$.\ The cases where\ $F\left|_{K}\right.\equiv 0$\ and where\ $K$\ is a singleton  are
trivial. Let\ $F\not\equiv 0$\ on\ $K$\ and let\ $[\alpha, \beta]$,\ $-\infty<\alpha<\beta<\infty$,\ be the closed
convex hull of\ $K$.\ Then Theorem 2.2 implies
\[
\max_{x\in K}|F(x)|=\max_{x\in [a,b]}|F(x)|=\max\{\ |F(\alpha), |F(\beta)|\}>0.
 \]
For definiteness, we suppose
\begin{equation}
\max_{x\in K}|F(x)|=\max_{x\in [\alpha, \beta]}|F(x)|=|F(\alpha)|>0.
 \end{equation}
Then
\begin{equation}
\varphi(x):=\frac{F(x+\alpha)}{F(\alpha)},\quad x\in [0, \beta-\alpha].
 \end{equation}
admits a positive-definite extension on\ $\R$.\ The \pd function\ $\varphi$\  satisfies\
$\varphi(-x)={\overline{\varphi(x)}}$,\ $x\in\R$,\ and coincides with continuous function\ $F(x+\alpha)/F(\alpha)$\ on\
$[-\beta-\alpha, \beta-\alpha]$.\ Thus\ (3.1)\ implies that\ $\varphi$\ is a continuous
 on the whole\ $\R$.\  Therefore\ $F\in CH(\R)$.\

\ \ \ \ \ \ Part 2. \ Let\ ${\cal A}$\ be a complex unital Banach algebra,\ $a\in H({\cal A})$,\ and\ $F\in\Hn$.\ As in
the proof of Part 1, we need only consider nontrivial cases where the spectrum\ $\sigma(a)$\  is not a singleton, and
where\ $F\left|_{\sigma(a)}\right.\not\equiv 0$.\ Similarly, we may without loss of generality  suppose that\ (3.10)\
also holds, where in this part\ $K=\sigma(a)$.\ It follows from Part 1     that the function\ (3.11)\ admits a
continuous positive-definite extension\ $\varphi$\  on\ $\R$.\ By Bochner's theorem, there is a positive measure\
$\eta\in\M$\ such that\ $\varphi(x)={\hat \eta}(x)$\ for all\ $x\in\R$.\  Thus\ (3.11)\ implies
\begin{equation}
F(x)=F(\alpha)\varphi(x-\alpha)=F(\alpha){\hat \eta}(x-\alpha),\quad  x\in [\alpha, \beta].
 \end{equation}
Now applying Lemma 3.2 to the symbols\ $F_1(x)=F(x)$,\ $F_2(x)=F(\alpha){\hat \eta}(x-\alpha)$,\ where\ $K$\ is replaced by the spectral synthesis set\
 $[\alpha,\beta]$,\ then gives\ $F(a)=F(\alpha){\hat \eta}(a-\alpha\mathbf{1}_{\cal A} )$.\ Combining\
 (1.2),\ (2.1)\ and\ (3.12),\ and noting that\ $\eta$\ is a
 probability measure, we  get
\begin{equation}
\|F(a)\|=\|F(\alpha)\varphi(a-\alpha\mathbf{1}_{\cal A})\|=|F(\alpha)|\|{\hat \eta}(a-\alpha\mathbf{1}_{\cal A})\|\le
|F(\alpha)|\|\eta\|_{\M}=|F(\alpha)|.
 \end{equation}
On the other hand,  it follows from the spectral mapping theorem\ (2.2)\  and\ (3.10)\ that the spectral radius of\
$F(a)$\ satisfies the equation
\[
|F(a)|=\max_{x\in \sigma(a)=K}|F(x)|=|F(\alpha)|.
\]
Now\ (3.13)\ implies that\ $\|F(a)\|=|F(a)|$.\ Therefore\ $F\in \U$.\ This completes the proof of Part 2 and Theorem
2.1.
\newpage

 \centerline{\bf{REFERENCES}}
 \vspace{4mm}

1. \  F.F. Bonsall and T. Duncan,\ {\it Numerical Ranges of Operators on Normed Spaces and of Elements of Normed
Algebras},\ Cambridge Univ. Press, London (1971).

2. \  F.F. Bonsall and T. Duncan,\ {\it Numerical Ranges II},\ Cambridge Univ. Press, London (1973).

3. \  E.A.  Gorin,\  Bernstein's inequality from the point of view of operator theory,\  Selecta Math. Sovietica, {\bf
7}(3), 91--219 (1988).

4. \  E.A.  Gorin,\  Universal symbols on locally compact abelian groups,\  Bull. Pol. Acad. Sci., Math.,  {\bf 51}(
2), 199--204 (2003).

5. \  E.A.  Gorin,\  Estimates for the involution of decomposable elements of a complex Banach algebra,\ Funct. Anal.
Appl., {\bf  39}( 4),  256--270 (2005); translation from Funkts. Anal. Prilozh., {\bf  39}(4), 14--31 (2005).

6. \ E. Hewitt and K.A. Ross, \ {\it Abstract Harmonic Analysis II},\ Springer-Verlag, Berlin-Heidelberg (1997).

7.\  H. K\"onig, \ A functional calculus for Hermitian elements of complex Banach algebras,\ Arch. Math., {\bf 28}(4),
422--430, (1977).

8. \  B.Ya.  Levin,\ {\it Lectures on Entire Functions},\ Translations of Mathematical Monographs, {\bf 150}, American
Math. Soc. (1996).

9. \ S. Norvidas,\ On stability of differential operators in spaces of entire functions,\  Sov. Math., Dokl.,  {\bf
34},  521-524 (1987); translation from Dokl. Akad. Nauk SSSR, {\bf 291}, 548--551 (1986).

10. \ S. Norvidas,\ Universal symbols of Hermitian elements,\ Preprint, {\bf 99-27},  1--12, Vilnius University, Faculty
Math. Inf., Vilnius (1999).

11. \ A.M. Sinclair,\ The norm of a Hermitian element in a Banach algebra,\ Proc. Amer. math. Soc., {\bf 28}(2),
446-450 (1971).

}}
\end{document}